\documentclass{amsart}

\input{epsf.tex}
\usepackage{amsfonts,amssymb,verbatim,amsmath,amsthm,latexsym,textcomp,amscd}
\usepackage{latexsym,amsfonts,amssymb,epsfig,verbatim}
\usepackage{amsmath,amsthm,amssymb,latexsym,graphics,textcomp}
\usepackage{graphicx}
\usepackage{color}
\usepackage{url}

\begin{document}

\newtheorem{theorem}{Theorem}[section]
\newtheorem{prop}[theorem]{Proposition}
\newtheorem{lemma}[theorem]{Lemma}
\newtheorem{cor}[theorem]{Corollary}
\newtheorem{defn}[theorem]{Definition}
\newtheorem{conj}[theorem]{Conjecture}
\newtheorem{claim}[theorem]{Claim}

\newcommand{\boundary}{\partial}
\newcommand{\bbC}{{\mathbb C}}
\newcommand{\bbD}{{\mathbb D}}
\newcommand{\bbH}{{\mathbb H}}
\newcommand{\bbZ}{{\mathbb Z}}
\newcommand{\bbN}{{\mathbb N}}
\newcommand{\bbQ}{{\mathbb Q}}
\newcommand{\bbR}{{\mathbb R}}
\newcommand{\proj}{{\mathbb P}}
\newcommand{\lhp}{{\mathbb L}}
\newcommand{\tube}{{\mathbb T}}
\newcommand{\cusp}{{\mathbb P}}
\newcommand\AAA{{\mathcal A}}
\newcommand\BB{{\mathcal B}}
\newcommand\CC{{\mathcal C}}
\newcommand\DD{{\mathcal D}}
\newcommand\EE{{\mathcal E}}
\newcommand\FF{{\mathcal F}}
\newcommand\GG{{\mathcal G}}
\newcommand\HH{{\mathcal H}}
\newcommand\II{{\mathcal I}}
\newcommand\JJ{{\mathcal J}}
\newcommand\KK{{\mathcal K}}
\newcommand\LL{{\mathcal L}}
\newcommand\MM{{\mathcal M}}
\newcommand\NN{{\mathcal N}}
\newcommand\OO{{\mathcal O}}
\newcommand\PP{{\mathcal P}}
\newcommand\QQ{{\mathcal Q}}
\newcommand\RR{{\mathcal R}}
\newcommand\SSS{{\mathcal S}}
\newcommand\TT{{\mathcal T}}
\newcommand\UU{{\mathcal U}}
\newcommand\VV{{\mathcal V}}
\newcommand\WW{{\mathcal W}}
\newcommand\XX{{\mathcal X}}
\newcommand\YY{{\mathcal Y}}
\newcommand\ZZ{{\mathcal Z}}
\newcommand\CH{{\CC\HH}}
\newcommand\TC{{\TT\CC}}
\newcommand\EXH{{ \EE (X, \HH )}}
\newcommand\GXH{{ \GG (X, \HH )}}
\newcommand\GYH{{ \GG (Y, \HH )}}
\newcommand\PEX{{\PP\EE  (X, \HH , \GG , \LL )}}
\newcommand\MF{{\MM\FF}}
\newcommand\PMF{{\PP\kern-2pt\MM\FF}}
\newcommand\ML{{\MM\LL}}
\newcommand\PML{{\PP\kern-2pt\MM\LL}}
\newcommand\GL{{\GG\LL}}
\newcommand\Pol{{\mathcal P}}
\newcommand\half{{\textstyle{\frac12}}}
\newcommand\Half{{\frac12}}
\newcommand\Mod{\operatorname{Mod}}
\newcommand\Area{\operatorname{Area}}
\newcommand\ep{\epsilon}
\newcommand\hhat{\widehat}
\newcommand\Proj{{\mathbf P}}
\newcommand\U{{\mathbf U}}
 \newcommand\Hyp{{\mathbf H}}
\newcommand\D{{\mathbf D}}
\newcommand\Z{{\mathbb Z}}
\newcommand\R{{\mathbb R}}
\newcommand\Q{{\mathbb Q}}
\newcommand\E{{\mathbb E}}
\newcommand\til{\widetilde}
\newcommand\length{\operatorname{length}}
\newcommand\tr{\operatorname{tr}}
\newcommand\gesim{\succ}
\newcommand\lesim{\prec}
\newcommand\simle{\lesim}
\newcommand\simge{\gesim}
\newcommand{\simmult}{\asymp}
\newcommand{\simadd}{\mathrel{\overset{\text{\tiny $+$}}{\sim}}}
\newcommand{\ssm}{\setminus}
\newcommand{\diam}{\operatorname{diam}}
\newcommand{\pair}[1]{\langle #1\rangle}
\newcommand{\T}{{\mathbf T}}
\newcommand{\inj}{\operatorname{inj}}
\newcommand{\pleat}{\operatorname{\mathbf{pleat}}}
\newcommand{\short}{\operatorname{\mathbf{short}}}
\newcommand{\vertices}{\operatorname{vert}}
\newcommand{\collar}{\operatorname{\mathbf{collar}}}
\newcommand{\bcollar}{\operatorname{\overline{\mathbf{collar}}}}
\newcommand{\I}{{\mathbf I}}
\newcommand{\tprec}{\prec_t}
\newcommand{\fprec}{\prec_f}
\newcommand{\bprec}{\prec_b}
\newcommand{\pprec}{\prec_p}
\newcommand{\ppreceq}{\preceq_p}
\newcommand{\sprec}{\prec_s}
\newcommand{\cpreceq}{\preceq_c}
\newcommand{\cprec}{\prec_c}
\newcommand{\topprec}{\prec_{\rm top}}
\newcommand{\Topprec}{\prec_{\rm TOP}}
\newcommand{\fsub}{\mathrel{\scriptstyle\searrow}}
\newcommand{\bsub}{\mathrel{\scriptstyle\swarrow}}
\newcommand{\fsubd}{\mathrel{{\scriptstyle\searrow}\kern-1ex^d\kern0.5ex}}
\newcommand{\bsubd}{\mathrel{{\scriptstyle\swarrow}\kern-1.6ex^d\kern0.8ex}}
\newcommand{\fsubeq}{\mathrel{\raise-.7ex\hbox{$\overset{\searrow}{=}$}}}
\newcommand{\bsubeq}{\mathrel{\raise-.7ex\hbox{$\overset{\swarrow}{=}$}}}
\newcommand{\tw}{\operatorname{tw}}
\newcommand{\base}{\operatorname{base}}
\newcommand{\trans}{\operatorname{trans}}
\newcommand{\rest}{|_}
\newcommand{\bbar}{\overline}
\newcommand{\UML}{\operatorname{\UU\MM\LL}}
\newcommand{\EL}{\mathcal{EL}}
\newcommand{\tsum}{\sideset{}{'}\sum}
\newcommand{\tsh}[1]{\left\{\kern-.9ex\left\{#1\right\}\kern-.9ex\right\}}
\newcommand{\Tsh}[2]{\tsh{#2}_{#1}}
\newcommand{\qeq}{\mathrel{\approx}}
\newcommand{\Qeq}[1]{\mathrel{\approx_{#1}}}
\newcommand{\qle}{\lesssim}
\newcommand{\Qle}[1]{\mathrel{\lesssim_{#1}}}
\newcommand{\simp}{\operatorname{simp}}
\newcommand{\vsucc}{\operatorname{succ}}
\newcommand{\vpred}{\operatorname{pred}}
\newcommand\fhalf[1]{\overrightarrow {#1}}
\newcommand\bhalf[1]{\overleftarrow {#1}}
\newcommand\sleft{_{\text{left}}}
\newcommand\sright{_{\text{right}}}
\newcommand\sbtop{_{\text{top}}}
\newcommand\sbot{_{\text{bot}}}
\newcommand\sll{_{\mathbf l}}
\newcommand\srr{_{\mathbf r}}
\newcommand\geod{\operatorname{\mathbf g}}
\newcommand\mtorus[1]{\boundary U(#1)}
\newcommand\A{\mathbf A}
\newcommand\Aleft[1]{\A\sleft(#1)}
\newcommand\Aright[1]{\A\sright(#1)}
\newcommand\Atop[1]{\A\sbtop(#1)}
\newcommand\Abot[1]{\A\sbot(#1)}
\newcommand\boundvert{{\boundary_{||}}}
\newcommand\storus[1]{U(#1)}
\newcommand\Momega{\omega_M}
\newcommand\nomega{\omega_\nu}
\newcommand\twist{\operatorname{tw}}
\newcommand\modl{M_\nu}
\newcommand\MT{{\mathbb T}}
\newcommand\Teich{{\mathcal T}}
\renewcommand{\Re}{\operatorname{Re}}
\renewcommand{\Im}{\operatorname{Im}}

\title{Simultaneous linearization of germs of commuting holomorphic diffeomorphisms}

\author{Kingshook Biswas }

\date{}

\thanks{Research partly supported by  Department of Science and Technology research project
grant DyNo. 100/IFD/8347/2008-2009}

\maketitle

\begin{abstract} Let $\alpha_1, \dots, \alpha_n$ be irrational numbers
which are linearly independent over the rationals, and $f_1,
\dots, f_n$ commuting germs of holomorphic diffeomorphisms in
$\mathbb{C}$ such that $f_k(0) = 0, f_k'(0) = e^{2\pi i \alpha_k},
k=1,\dots,n$. Moser showed that $f_1, \dots, f_n$ are
simultaneously linearizable (i.e. conjugate by a germ of
holomorphic diffeomorphism $h$ to the rigid rotations
$R_{\alpha_k}(z) = e^{2\pi i \alpha_k}z$) if the vector of
rotation numbers $(\alpha_1, \dots, \alpha_n)$ satisfies a
Diophantine condition. Adapting Yoccoz' renormalization to
the setting of commuting germs, we show that simultaneous linearization
holds in the presence of a weaker Brjuno-type condition ${\mathcal
B}(\alpha_1, \dots, \alpha_n) < +\infty$ where ${\mathcal
B}(\alpha_1, \dots, \alpha_n)$ is a multivariable analogue of the
Brjuno function. If there are no periodic orbits for the
action of the germs $f_1, \dots, f_n$ in a neighbourhood of the
origin then a weaker arithmetic
condition ${\mathcal B}'(\alpha_1, \dots, \alpha_n) < +\infty$ analogous
to Perez-Marco's condition for linearization in the absence of periodic
orbits is shown to suffice for linearizability. Normalizing the germs
to be univalent on the unit disk, in both cases the
Siegel disks are shown to contain disks of radii $C e^{-2\pi
\mathcal{B}}, C e^{-2\pi \mathcal{B'}}$ respectively for some
universal constant $C$.

\smallskip

\begin{center}

{\em AMS Subject Classification: 37F50}

\end{center}

\end{abstract}

\overfullrule=0pt

\tableofcontents


\section{Introduction}

\medskip

A germ $f(z) = e^{2\pi i \alpha}z + O(z^2), \alpha \in \mathbb{R - Q}/\mathbb{Z}$
of holomorphic diffeomorphism fixing the origin in $\mathbb{C}$ is
said to be linearizable if it is analytically conjugate to the
rigid rotation $R_{\alpha}(z) = e^{2\pi i \alpha} z$. The number $\alpha$
is called the rotation number of $f$, and the maximal domain of linearization
is called the Siegel disk of $f$. Siegel \cite{siegel}
showed that if $\alpha$ satisfies a Diophantine condition
$$
\left|\alpha - \frac{p}{q}\right| \geq \frac{C}{q^\tau}
$$
for all integers $p,q$ for some $C, \tau >0$, then $f$ is
linearizable. Brjuno \cite{brjuno} showed that linearizability
holds under the more general Brjuno condition defined as follows:
Let $(\alpha^{(i)})_{i \geq 0}$ denote the fractions arising from the
continued fraction algorithm (given by iterating the Gauss map
$\alpha \mapsto 1/\alpha (mod 1)$) applied to $\alpha$,
then the Brjuno condition is
$$
\mathcal{B}(\alpha) := \sum_{n = 0}^{\infty} {\alpha}^{(0)} \dots
{\alpha}^{(n-1)} \log \left(\frac{1}{{\alpha}^{(n)}}\right) < +\infty.
$$
Yoccoz \cite{yoccoz} 
proved the optimality of the Brjuno condition: if $\mathcal{B}(\alpha) = +\infty$
then there exists a nonlinearizable germ $f$ with rotation number
$\alpha$. Perez-Marco \cite{perezmens} showed that the weaker arithmetic condition
$$
\mathcal{B}'(\alpha) := \sum_{n = 0}^{\infty} {\alpha}^{(0)} \dots
{\alpha}^{(n-1)} \log \log \left(\frac{e}{{\alpha}^{(n)}}\right) < +\infty.
$$
is sufficient for linearizability in the absence of periodic orbits accumulating the fixed
point, and is optimal in this case: if $\mathcal{B}'(\alpha) = +\infty$
then there exists a nonlinearizable germ $f$ with rotation number
$\alpha$ with no periodic orbits in a neighbourhood of the origin.

\medskip

It is well known that $N$ commuting germs $f_1, \dots, f_N$ are
either simultaneously linearizable or nonlinearizable. Moser
\cite{moser} has shown that if the rotation numbers $\alpha_1, \dots, \alpha_N$
satisfy a Diophantine condition
$$
\max_{1 \leq k \leq N} |q \alpha - p_k| \geq \frac{C}{q^{\tau}}
$$
for some $C>0, \tau>0$ then the germs $f_1, \dots, f_N$ are simultaneously
linearizable. More precisely Moser proves that commuting analytic circle diffeomorphisms whose
rotation numbers satisfy the above condition and are sufficiently close to
the corresponding rotations are simultaneously
analytically linearizable; by the dictionary between circle maps and germs established
by Perez-Marco in \cite{perezmcirclemaps} 
this local linearization result implies the corresponding global linearization
result for germs.

\medskip

Generalizing Moser's result we show that multidimensional analogues of Brjuno's and Perez-Marco's conditions
given by arithmetic functions $\mathcal{B}(\alpha_1,\dots,\alpha_N)$,
$\mathcal{B}'(\alpha_1,\dots,\alpha_N)$ (defined in the next section) suffice for
linearizability in the presence and absence of periodic
orbits:

\medskip

\begin{theorem} \label{linearization1}{Let $f_1, \dots, f_N$
be $N$ commuting germs with rationally independent irrational rotation numbers
$\alpha_1,\dots,\alpha_N$ univalent on the unit disk $\mathbb{D}$.
\ If $\mathcal{B}(\alpha_1, \dots, \alpha_N) < +\infty$ then
$f_1,\dots,f_N$ are simultaneously linearizable, and their common Siegel
disk contains the disk of radius $C e^{-2\pi \mathcal{B}(\alpha_1, \dots, \alpha_N)}$
around the origin for some universal constant $C$.}
\end{theorem}

\medskip

In the appendix it is shown that the above arithmetic condition is implied
by Moser's Diophantine condition.
By a periodic orbit for the action of collection of commuting maps $f_1,\dots,
f_N$ univalent on $\mathbb{D}$ we mean a finite set $\mathcal{O} \subset
\mathbb{D}$ invariant under the maps $f_j$, i.e. $f_j(\mathcal{O}) = \mathcal{O},
j=1,\dots,N$.

\medskip

\begin{theorem} \label{linearization2}{Let $f_1, \dots, f_N$
be $N$ commuting germs with rationally independent irrational rotation numbers
$\alpha_1,\dots,\alpha_N$, univalent on the unit disk $\mathbb{D}$
with no periodic orbits in $\mathbb{D}$.
If $\mathcal{B}'(\alpha_1, \dots, \alpha_N) < +\infty$ then
$f_1,\dots,f_N$ are simultaneously linearizable, and their common Siegel
disk contains the disk of radius
$C e^{-2\pi \mathcal{B}'(\alpha_1, \dots, \alpha_N)}$ around the origin for some
universal constant $C$.}
\end{theorem}

\medskip

%
%


\section{Brjuno-type functions.}

\medskip

We define Gauss maps for $N$-tuples of fractions as follows:

\medskip

Let $\alpha_1, \dots, \alpha_N \in (0,1)$ be irrationals linearly
independent over $\bbQ$, and $w \in \{1,\dots,N\}$. We define
irrationals $(\tilde{\alpha_1},\dots,\tilde{\alpha_N}) = G(\alpha_1,\dots,\alpha_N,w),
0< \tilde{\alpha_i} < 1/2$, as follows:

\medskip

Let $\hat{\alpha_i} = \alpha_i$ for $i \neq w$ and $\hat{\alpha_i}
= 1$ for $i = w$.

\medskip

\noindent Let $a_i = a_i(\alpha_1,\dots, \alpha_N,w)$ be integers such that

$$
|a_i \alpha_w - \hat{\alpha}_i| = \min_{k \in \bbZ} |k\alpha_w - \hat{\alpha}_i| < \frac{1}{2}\alpha_w
$$

\noindent Let
$$
\tilde{\alpha}_i := \frac{|a_i \alpha_w - \hat{\alpha}_i|}{\alpha_w} \in (0, 1/2),
\epsilon_i = \hbox{sgn}(a_i \alpha_w - \hat{\alpha}_i)
$$

\medskip

Given an infinite word $\overline{w} = (w(n))_{n \in \mathbb{N}} \in
\mathcal{C} := \{1,\dots,N\}^{\mathbb{N}}$, let $\alpha^{(0)}_i
= \alpha_i$ and
$(\alpha^{(n+1)}_1,\dots,\alpha^{(n+1)}_N) = G(\alpha^{(n)}_1,\dots,\alpha^{(n)}_N,
w(n))$ for $n \geq 0$. We define
\begin{align*}
\mathcal{B}(\alpha_1, \dots, \alpha_n, \overline{w}) & :=  \sum_{n = 0}^{\infty}
\alpha^{(0)}_{w(0)}\dots\alpha^{(n-1)}_{w(n-1)} \log (\frac{1}{\alpha^{(n)}_{w(n)}}) \\
\mathcal{B'}(\alpha_1, \dots, \alpha_n, \overline{w}) & :=  \sum_{n = 0}^{\infty}
\alpha^{(0)}_{w(0)}\dots\alpha^{(n-1)}_{w(n-1)} \log \log (\frac{1}{
\alpha^{(n)}_{w(n)}}) \\
\end{align*}
and
\begin{align*}
\mathcal{B}(\alpha_1, \dots, \alpha_n) & := \inf_{\overline{w} \in \mathcal{C}} \mathcal{B}(\alpha_1, \dots, \alpha_n,
\overline{w}) \\
\mathcal{B'}(\alpha_1, \dots, \alpha_n) & := \inf_{\overline{w} \in \mathcal{C}} \mathcal{B'}(\alpha_1, \dots, \alpha_n,
\overline{w}) \\
\end{align*}

\bigskip

\section{Renormalization of commuting germs.}

\medskip

In this section we generalize Yoccoz' renormalization for germs 
to the setting of $N$ commuting germs. Throughout $C>0$ will denote various universal
constants, ${\bbD}_r$ the disc $\{ |z| < r \}$, $\mathbb{H}_y$ the
half-plane $\{ \Im z > y \}$ and $S(a,y)$ the vertical strip $\{ |\Re z|
\leq a, \Im z \geq y \}$.

\medskip

\medskip

Let
$$
\mathcal{S}_{t,r} := \{ \ f \ | \ f \hbox{ is univalent on } \bbD, f(0) = 0,
f'(0) = e^{2\pi i \alpha} \ \}
$$
and
$$
\hat{\mathcal{S}}_{t,y} := \{ \ F \ | \ F \hbox{ is univalent on } \mathbb{H}_y,
\lim_{\Im z \to +\infty} (F(z) - (z + \alpha)) = 0, F(z+1) = F(z)+1 \ \}
$$
be the corresponding space of lifts to the half-plane $\mathbb{H}_y$ via the
universal covering $E : \mathbb{H} \to \mathbb{D}^*, E(z) =
e^{2\pi i z}$, where $y = \frac{1}{2\pi}\log(1/r)$. These spaces are
compact (with respect to the topology of uniform
convergence on compacta), and satisfy
\begin{align*}
|F(z) - (z+\alpha)| & \leq C e^{-2\pi \Im z} \\
|F'(z) - 1| & \leq C e^{-2\pi \Im z} \\
\end{align*}
for all $F \in \hat{\mathcal{S}}_{\alpha}, \Im z \geq C + y$. We call $\alpha$ the rotation
number of $F$. By the Koebe one-quarter theorem,
the inverse $F^{-1}$ is well-defined and univalent in $\mathbb{H}_{C+y}$
for a universal constant $C$, and satisfies the same estimates
as $F$ for $C$ large enough.

\medskip

\subsection{Uniformizing maps}

\medskip

We recall the uniformizing maps defined by Yoccoz in order to
perform the renormalization. 
For $\alpha \in \mathbb{R}$ let $t(\alpha) > 0$ be such
that $|F(z) - (z+\alpha)| \leq \frac{\alpha}{4}$ for $\Im z \geq
t(\alpha), F \in \hat{\mathcal{S}}_{\alpha,0}$ (by the above
estimates we can take $t(\alpha) = \frac{1}{2\pi}\log(\frac{1}{\alpha}) + C$
for a universal constant $C$).

\medskip

Let $l$ be the line segment $[it(\alpha), +i\infty]$. Then $l, F(l)$ and the segment
$s = [it(\alpha), F(it(\alpha))]$ bound a domain $\mathcal{U}$ in
$\mathbb{H}$. The proximity of $F$ to the translation by
$\alpha$ shows that any point $z$ in the domain
$\mathcal{V} = \{ \Im z \geq \frac{1}{2}|\Re z| + t(\alpha) + 1 \}$ has an
iterate $F^n(z)$ contained in $\mathcal{U}$.

\medskip

Pasting the boundary segments $l$ and $F(l)$ of $\mathcal{U}$ by
$F$ the quotient $\pi : \overline{\mathcal{U}} \to \mathcal{S} := \overline{\mathcal{U}}/(l \sim F(l))$
gives a Riemann surface $\mathcal{S}$ which is biholomorphic
to the punctured disc $\mathbb{D}^*$. Let $k : \mathcal{S} \to \mathbb{D}^*$
be a uniformization and $K : \overline{\mathcal{U}} \to \mathbb{H}$ a lift
of $k$ satisfying $E \circ K = k \circ \pi$. The lift $K$
satisfies $K(F(z)) = K(z) + 1$ for points $z$ in the left boundary
segment $l$ of $\mathcal{U}$. We can therefore extend $K$ univalently to
$\mathcal{V}$ by defining $K(z) := K(F^j(z)) - j$ where $z \in \mathcal{V},
F^j(z) \in \mathcal{U}$. The extended map $K$ satisfies
$K(F(z)) = K(z) + 1$ whenever $z, F(z) \in \mathcal{V}$ and is called a {\it
uniformizing map} for $F$.

\medskip

For $C > 0$ large enough the domain $\mathcal{U'}$ bounded by
$l'=[i(t(\alpha)+C), +i\infty], F(l')$ and $s' = [i(t(\alpha)+C), F(i(t(\alpha)+C))]$
is contained in $\mathcal{V}$ and Yoccoz shows that 
the following estimate holds:

\medskip

\begin{prop}\label{derivative close to constant}{For
$z \in S(2, t(\alpha)+C) \subset \mathcal{V}$
$$
\left|K'(z) - \frac{1}{\alpha}\right|  \leq \frac{1}{3}\frac{1}{(y - t(\alpha))^2}
$$
}
\end{prop}

\medskip

Composing $K$ with a translation we may assume that $K(s')$ is
contained in the closure of the lower half-plane with
$K(s') \cap \mathbb{R} \neq \emptyset$. The above estimate then
yields (for $C$ large enough) on integration

\medskip

\begin{prop}\label{close to affine}{For
$z \in S(2, t(\alpha)+C)$,
$$
\left|K(z) - \frac{1}{\alpha}\left(z - i(t(\alpha)+C)\right)\right| \leq
\frac{2}{\alpha}
$$
}
\end{prop}

\medskip

We may assume $C$ is large enough so that for $|j| \leq
a(\alpha) = [10/\alpha]+1$, the maps $F^j$ are univalent on
$\mathbb{H}_{t(\alpha) + C}$ mapping into $\mathbb{H}_{t(\alpha) + C -
1}$. For $n \leq a(\alpha)$ let
$$
\mathcal{A}(F, n) = \bigcup_{j = -n}^{n} F^j(\mathcal{U}')
$$
Taking $C$ large enough, $\mathcal{A}(F,a(\alpha))$ contains the strip
$S(8, t(\alpha) + C + 1)$.
The map $K$ extends to $\mathcal{A}(F, n)$ putting $K(F^j(z)) =
K(z) + j$, and by the above proposition the image $K(l')$ of the vertical line
$l'$ is close enough to the positive imaginary axis so that
$$
K(\mathcal{A}(F,a(\alpha)) = \bigcup_{j = -a(\alpha)}^{j = a(\alpha)}
T^j(K(\mathcal{U}'))
$$
covers the strip $S(8/\alpha,0)$.

\medskip

\subsection{Simultaneous renormalization.}

\medskip

Given $N$ commuting maps $F_j \in \hat{\mathcal{S}}_{\alpha_j}, j = 1,\dots,N$
with irrational rationally independent rotation
numbers $\alpha_j, j=1,\dots,N$ and $w \in \{1,\dots,N\}$ we
define a renormalization operator $\mathcal{R}$ giving $N$
commuting maps $(\tilde{F_1}, \dots, \tilde{F_N}) = \mathcal{R}(F_1,\dots,F_N,w)$ with rotation numbers
$(\tilde{\alpha}_1,\dots,\tilde{\alpha}_N) = G(\alpha_1,\dots,\alpha_N,w)$ as
follows:

\medskip

Let $K$ be a uniformizing map for $F = F_w$ with associated domains
$\mathcal{U}, \mathcal{U'}, \mathcal{V}$ as defined above.

\medskip

Let $\hat{F}_i = F_i$ for $i \neq w$, and $\hat{F}_i = T$ for $i =
w$, where $T(z) = z+1$ is the translation by one.

\medskip

We may assume the universal constant $C>2$ is
chosen large enough so that the maps $F^{\epsilon_i a_i}_w \circ
\hat{F}^{-\epsilon_i}_i$ are univalent on
the strip $S(2, t(\alpha)+C)$ and map it into $\mathcal{V}$.
On the domain $\mathcal{X} := K(S(2, t(\alpha)+C))$ we define
the maps
$$
\tilde{F}_i := K \circ (F^{\epsilon_i a_i}_w \circ
\hat{F}^{-\epsilon_i}_i) \circ K^{-1}
$$
Since $K$ conjugates $F_w$ to the translation $T$, the maps $\tilde{F}_i$
commute with $T$ and therefore extend univalently to the union
of translates $\mathcal{Y} := \cup_{n \in \mathbb{Z}} T^n(\mathcal{X})$
by setting $\tilde{F}_i(T^n(z)) := \tilde{F}_i(z) + n$ for $z \in \mathcal{X}$. Since the image $K(s')$
of the segment $s'$ lies in the lower half-plane, $\mathcal{Y}$ contains the upper half-plane,
so we obtain $N$ commuting maps $\tilde{F}_j \in \hat{\mathcal{S}}_{\tilde{\alpha}_j, 0}, j=1,\dots,N$.

\medskip

\begin{lemma}\label{pulling back}{Given $G \in \hat{\mathcal{S}}_{\tilde{\beta},\tilde{y}}$
commuting with $\tilde{F}_1, \dots, \tilde{F}_N$, where $|\tilde{\beta}| < 1/2, \tilde{y} \geq 0$,
the conjugate $G^* = K^{-1} \circ G \circ K$ defines an element of
$\hat{\mathcal{S}}_{\beta,y}$ commuting with $F_1,\dots,F_N$ where
$\beta = \alpha_w \tilde{\beta}$ and $y = \alpha_w \tilde{y} + t(\alpha_w) + (1 + \alpha_w)C + 2$.
}
\end{lemma}

\medskip

\noindent{\bf Proof:} It follows from Proposition \ref{close to affine} that
the image of the strip $S(2/\alpha_w, y)$
under $K$ is contained in the strip $S(4/\alpha_w,\tilde{y}+C)$. For
$C$ large enough the maps $G, \tilde{F}^{\epsilon}_j \circ T^k, G \circ
\tilde{F}^{\epsilon}_j \circ T^k,
\epsilon= \pm 1, j=1\,\dots,N, |k| \leq [1/\alpha_w]+1,$ map
$S(4/\alpha_w,\tilde{y}+C)$ into $S(8/\alpha_w,\tilde{y}) \subset
K(\mathcal{A}(F_w, a(\alpha_w))$. So the conjugates of these maps
by $K^{-1}$ are univalent in $S(2/\alpha_w, y)$; the conjugates
include $G^*, F_1, \dots, F_N, T$ and commute because
$G \circ
\tilde{F}^{\epsilon}_j \circ T^k = \tilde{F}^{\epsilon}_j \circ T^k \circ
G$. Since $G^*$ commutes with $T$, it extends to a univalent map
on the union of translates of $S(2/\alpha_w,y)$ under $T$ which contains
$\mathbb{H}_{y}$. From Proposition \ref{close to affine}, $K(z) =
z/\alpha_w + o(z), K^{-1}(z) = \alpha_w z + o(z)$ as $\Im z \to
+\infty$ (with $\Re z$ bounded), from which it follows that
$$
\lim_{\Im z \to +\infty} G^*(z) - z = \alpha_w \tilde{\beta}
$$
so $G^* \in \hat{\mathcal{S}}_{\beta,y}$. $\diamond$

\medskip

\subsection{Iterating the renormalization}

\medskip

Given maps $F_1,\dots,F_N$ as above and an infinite
word $\overline{w} = (w(n))_{n \geq 0} \in \{1,\dots,N\}^{\mathbb{N}}$, we define
\begin{align*}
(F^{(0)}_1,\dots,F^{(0)}_N) & := (F_1,\dots,F_N) \\
(F^{(n+1)}_1,\dots,F^{(n+1)}_N) & := \mathcal{R}(F^{(n)}_1,\dots,F^{(n)}_N, w(n)) \ , \ n \geq 0 \\
\end{align*}

We denote by $\tilde{F}_n$ the map $F^{(n)}_{w(n)}$
with respect to which we
are renormalizing at the $n$th stage and by $\tilde{\alpha}_n$ its
rotation number $\alpha^{(n)}_{w(n)}$. Let $K_n$ be the
uniformizing map of $\tilde{F}_n$

\medskip

For $0 \leq m \leq n$, we define finite families of commuting maps
$\mathcal{F}_{m,n}$ univalent on $\mathbb{H}_{y_{m,n}}$ as follows:

\medskip

\noindent Let $y_{n,n} := t(\tilde{\alpha}_n)$. Let $\mathcal{F}_{n,n}$ be the
collection of maps $\{ \tilde{F}^j_n : 0 \leq j \leq [1/\tilde{\alpha}_n] \}$ univalent
on $\mathbb{H}_{y_{m,n}}$. Let $\epsilon_{n,n} :=
\tilde{\alpha}_n$.

\medskip

\noindent For $0 < m \leq n$ assume $y_{m,n}$ and $\mathcal{F}_{m,n}$
have been defined and satisfy the following induction hypotheses:

\medskip

\noindent (1) The maps in $\mathcal{F}_{m,n}$ are univalent in
$\mathbb{H}_{y_{m,n}}$.

\noindent (2) $F^{(m)}_1, \dots, F^{(m)}_N, T$ commute
with the elements of $\mathcal{F}_{m,n}$.

\noindent (3) The rotation numbers of the maps in
$\mathcal{F}_{m,n}$ are contained in $[0, 1)$ and are
$\epsilon_{m,n}$-dense in $[0,1]$ (every point in $[0,1]$ is at
distance less than $\epsilon_{m,n}$ from these rotation numbers).

\medskip

Let $y_{m-1,n} := \tilde{\alpha}_{m-1} y_{m,n} + t(\tilde{\alpha}_{m-1}) + (1 +
\tilde{\alpha}_{m-1})C + 2$. By Lemma \ref{pulling back}, the collection of conjugates
$\tilde{\mathcal{F}}_{m-1,n} := \{ G^* = K^{-1}_{m-1} \circ G \circ K_{m-1} : G \in \mathcal{F}_{m,n} \}$
are univalent on $\mathbb{H}_{y_{m-1,n}}$, commute with $F^{(m-1)}_1, \dots, F^{(m-1)}_N, T$
and taking $C$ large enough we may assume
(since they have rotation numbers in $[0,1)$) that they map this half-plane into
$\mathbb{H}_{y_{m-1,n}-2}$. For $0 \leq j \leq
[1/\tilde{\alpha}_{m-1}]$, the maps $\tilde{F}^j_{m-1}$ are univalent on
 $\mathbb{H}_{\tilde{\alpha}_{m-1} + C} \supset \mathbb{H}_{y_{m-1,n}-2}$;
 for each $G \in \tilde{\mathcal{F}}_{m-1,n}$ we can choose $0 \leq j(G^*) \leq [1/\tilde{\alpha}_{m-1}]$
such that the rotation numbers of the maps $\{ \tilde{F}^j_{m-1}
\circ G^* : G^* \in \tilde{\mathcal{F}}_{m-1,n}, 0 \leq j \leq
j(G^*) \}$ are contained in $[0,1)$, are $\tilde{\alpha}_{m-1}
\epsilon_{m,n}$-dense in $[0,1]$, and these maps are univalent on $\mathbb{H}_{y_{m-1,n}}$.

\medskip

Let $\mathcal{F}_{m-1,n} = \{ \tilde{F}^j_{m-1}
\circ G^* : G^* \in \tilde{\mathcal{F}}_{m-1,n}, 0 \leq j \leq
j(G^*) \}$ be this collection of maps, and
$\epsilon_{m-1,n} = \tilde{\alpha}_{m-1}
\epsilon_{m,n}$. Then the induction hypotheses (1), (2), and (3)
above are satisfied for $y_{m-1,n}, \mathcal{F}_{m-1,n}$ and
$\epsilon_{m-1,n}$. By descending induction on $m$ we get:

\medskip

\begin{prop}\label{heights}{For any word $(w(n))_{n \geq 0} \in \{1,\dots,N\}^{\mathbb{N}}$
there are families of maps $\mathcal{F}_{0,n}$ univalent on
$\mathbb{H}_{y_{0,n}}$ commuting with $F_1,\dots, F_N, T$, where
$$
y_{0,n} \leq \sum_{j = 0}^{n} \alpha^{(0)}_{w(0)} \dots
\alpha^{(j-1)}_{w(j-1)}
t\left({\alpha^{(j)}_{w(j)}}\right) + C'
$$
for a universal constant $C'$. Moreover the rotation numbers of the maps $F \in \mathcal{F}_{0,n}$
are $\epsilon_{0,n}$-dense in $[0,1]$ where
$\epsilon_{0,n} = {\alpha}^{(0)}_{w(0)} \dots
{\alpha}^{(n)}_{w(n)}$.}
\end{prop}

\medskip

\noindent{Proof:} Assume that $y_{m+1,n}$ satisfies
$$
y_{m+1,n} \leq \sum_{j = m+1}^{n} \alpha^{(m+1)}_{w(m+1)} \dots
\alpha^{(j-1)}_{w(j-1)}
t\left({\alpha^{(j)}_{w(j)}}\right) + 2(C+1)\left(1 +
\sum_{j=m+1}^{n-2} \alpha^{(j)}_{w(j)} \dots \alpha^{(n-2)}_{w(n-2)}\right)
$$
(for an $m+1 \leq n$). Then the definition of $y_{m,n}$ gives

\begin{align*}
y_{m,n} \leq & \sum_{j = m}^{n} \alpha^{(m)}_{w(m)} \dots
\alpha^{(j-1)}_{w(j-1)}
t\left({\alpha^{(j)}_{w(j)}}\right) + \alpha^{(m)}_{w(m)} 2(C+1)\left(1 +
\sum_{j=m+1}^{n-2} \alpha^{(j)}_{w(j)} \dots
\alpha^{(n-2)}_{w(n-2)}\right) \\
             & + \left(1+\alpha^{(m)}_{w(m)}\right)C + 2 \\
        \leq & \sum_{j = m}^{n} \alpha^{(m)}_{w(m)} \dots
\alpha^{(j-1)}_{w(j-1)}
t\left({\alpha^{(j)}_{w(j)}}\right) + 2(C+1)\left(1 +
\sum_{j=m}^{n-2} \alpha^{(j)}_{w(j)} \dots
\alpha^{(n-2)}_{w(n-2)}\right) \\
        & \leq \sum_{j = m}^{n} \alpha^{(m)}_{w(m)} \dots
\alpha^{(j-1)}_{w(j-1)}
t\left({\alpha^{(j)}_{w(j)}}\right) + C'
\end{align*}
where $C' = 4(C+1)$, using $\alpha^{(j)}_{w(j)} < 1/2$
for all $j$. By induction we get the required inequality for $y_{0,n}$. $\diamond$

\medskip

\section{Linearization of commuting germs}

\medskip

\begin{lemma}\label{flow}{Let $\mathcal{F} = (f_t)_{t \in A}$ be a commuting family
of maps univalent on a disc $\mathbb{D}_r$ such that $f_t(0) = 0,
f'_t(0) = e^{2\pi i t}$, the set of rotation numbers $A$ contains
an irrational and is dense in $[0,1]$. Then the maps $f_t$ are simultaneously
linearizable with Siegel disk containing $\mathbb{D}_r$.}
\end{lemma}

\medskip

\noindent{Proof:} The family $\mathcal{F}$, being a normalized family of univalent functions
on a disk, is normal. Fix a $t_0$ in $A$ which is irrational.
For any sequence $t_n \in A$ converging to $t \in [0,1]$,
any normal limit $f$ of the sequence $(f_{t_n})$ must have
rotation number $t$ and commute with $f_{t_0}$; but $t_0$ being irrational,
the power series of $f$ is uniquely determined by these two conditions. So the sequence
$(f_{t_n})$ has a unique normal limit and is normally convergent to a map which we
denote by $f_t$. Since $A$ is dense we get a commuting family of maps $\{ f_t : t \in [0,1]
\}$ univalent on $\mathbb{D}_r$; as above, $t_n \to t$ implies
$f_{t_n}(z) \to f_t(z)$ for $t_n,t \in [0,1]$ so $t \mapsto
f_t(z)$ is continuous for fixed $z$. Let
$$
h(z) := \int_{0}^{1} e^{-2\pi i t} f_t(z) \ dt \ , \ z \in
\mathbb{D}_r.
$$
Then $h$ is holomorphic, $h(0) = 0$ and $h'(0) = 1$ (Riemann sums for the
integral converge normally to $h$).
For $t \in \mathbb{R}$ we let $f_t := f_{\overline{t}}$ where $\overline{t} \in [0,1]$
is the fractional part of $t$. Then the maps $f_t$ satisfy $f_t \circ f_s =
f_{t+s}$ for all real $t,s$, which gives $h(f_t(z)) = e^{2\pi i t}
f_t(z)$ so $h$ linearizes the maps $f_t$. $\diamond$

\medskip

\medskip

\noindent{Proof of Theorem \ref{linearization1}}: Let $F_1, \dots,
F_N$ be lifts of $f_1, \dots, f_N$ with rotation numbers
$0 < \alpha_1,\dots, \alpha_N < 1$. Let $\overline{w} = (w(n))$ be a word such
$\mathcal{B}(\alpha_1,\dots,\alpha_N,\overline{w}) \leq \mathcal{B}(\alpha_1,\dots,\alpha_N) + 1$.
 Since $t(\alpha^{(n)}_{w(n)}) = \frac{1}{2\pi}\log(1/\alpha)+C$ and $0 < \alpha^{(n)}_{w(n)} < 1/2$
for $n \geq 1$,
$$
\sum_{n = 0}^{\infty} \alpha^{(0)}_{w(0)} \dots
\alpha^{(n-1)}_{w(n-1)}
t\left({\alpha^{(j)}_{w(j)}}\right) \leq
\mathcal{B}(\alpha_1,\dots,\alpha_N,\overline{w}) + C''
$$
for a universal constant $C''$. Let $\mathcal{F}_{0,n}$ be the
families of maps given by Proposition \ref{heights} and
$\mathcal{F} = \{ f : E \circ F = f \circ E \ \hbox{for some } F \in \mathcal{F}_{0,n}, n \geq 0 \
\}$. Then by the estimate in Proposition \ref{heights}, the family
$\mathcal{F}$ satisfies the hypotheses of Lemma \ref{flow} on a
disc $\mathbb{D}_r$ where $r = C''' e^{-2\pi
\mathcal{B}(\alpha_1,\dots,\alpha_N)}$ for some universal constant
$C'''$, these maps commute with
$f_1,\dots,f_N$, hence $f_1,\dots,f_N$ are linearizable with
common Siegel disk containing $\mathbb{D}_r$. $\diamond$

\medskip

\medskip

\noindent{Proof of Theorem \ref{linearization2}}: For a map $F \in
\hat{\mathcal{S}}_{\alpha, 0}$ without fixed points, as Perez-Marco shows in
\cite{perezmens}, we can take $t(\alpha) = \frac{1}{2\pi}\log \log
(e/\alpha)$ as the height above which $F$ is close to translation
by $\alpha$. We note that the set of fixed points in $\mathbb{H}_C$ (where $C > 0$)
of any $F_j$ is finite modulo the translation $T$, and hence, if
nonempty, gives rise to a periodic orbit (modulo $T$) for
$F_1,\dots,F_N,T$. Since $F^{(n+1)}_{1}, \dots, F^{(n+1)}_{N}, T$ are conjugates by $K_n$ of the group
generated by $F^{(n)}_1,\dots,F^{(n)}_N, T$, a periodic orbit of $F^{(n+1)}_{1}, \dots, F^{(n+1)}_{N}, T$
corresponds to a periodic orbit for $F^{(n)}_{1}, \dots, F^{(n)}_{N}, T$.
It follows that if $F_1, \dots, F_N, T$ have no periodic
orbits then neither do $F^{(n)}_1,\dots,F^{(n)}_N, T$.  in
particular $F^{(n)}_{w(n)}$ has no fixed points modulo $T$, and so
we can take $t(\alpha^{(n)}_{w(n)}) = \frac{1}{2\pi}\log \log
(e/\alpha^{(n)}_{w(n)})$ for all $n$. The proof then follows the
same lines as that of Theorem \ref{linearization1} above.
$\diamond$

\medskip

\section{Appendix}

\medskip

We prove that Moser's Diophantine condition implies the
Brjuno-type condition:

\medskip

Let $N \geq 2$ be an integer and $\tau > 0$. We will write $a \ll b$ for quantities
satisfying $a \leq k b$ where $k$ is a constant depending only on $N, \tau$, and
$a \asymp b$ for $a \ll b, b \ll a$.
For $C > 0$ let $DC_N(C, \tau)$ be the set of $N$-tuples
$(\alpha_1,\dots,\alpha_N)$ such that $\alpha_1,\dots,\alpha_N$
are rationally independent irrationals in the interval $(0,1)$
satisfying
$$
\max_{1 \leq j \leq N} |q \alpha_j - p_j| \geq
\frac{C}{|q|^{\tau}}
$$
for all integers $p_1, \dots, p_N, q$ with $q \neq 0$. Using Khintchine's Transference
Principle it can be shown (see \cite{diophapp}, Ch.IV)
that any $(\alpha_1,\dots,\alpha_N) \in DC_N(C, \tau)$ satisfies
$$
|(\alpha_1,\dots,\alpha_N,1) \cdot (p_1,\dots,p_N,q) | \geq \frac{C'}{||(p_1,\dots,p_N,q)||^{\tau'}}
$$
for all $(p_1, \dots, p_N, q) \in \mathbb{Z}^{N+1} - \{0\}$, where $||\cdot||$ denotes
the sup norm on $\mathbb{R}^{N+1}$, $C' \gg C$ and $\tau' = N\tau + N - 1$. Conversely if an inequality
$$
|(\alpha_1,\dots,\alpha_N,1) \cdot (p_1,\dots,p_N,q) | \geq \frac{C'}{||(p_1,\dots,p_N,q)||^{\tau'}}
$$
holds for all $(p_1, \dots, p_N, q) \in \mathbb{Z}^{N+1} - \{0\}$
then $(\alpha_1,\dots,\alpha_N) \in DC_N(C, \tau)$ where $C
\gg C'$ and $\tau = \frac{\tau' + 1 - N}{N}$.

\medskip

Letting $(\alpha_1,\dots,\alpha_N) \mapsto
(\tilde{\alpha_1},\dots,\tilde{\alpha_N}) = G(\alpha_1,\dots,\alpha_N,
w)$ denote the Gauss maps defined in section 2 (where $w \in \{1,\dots,N\}$), we have

\medskip

\begin{lemma}\label{dioph}{
$$
(\alpha_1,\dots,\alpha_N) \in DC_N(C, \tau) \Rightarrow
(\tilde{\alpha_1},\dots,\tilde{\alpha_N}) \in DC_N(C', \tau)
$$
where $C' \gg C \alpha^{\tau' - 1}_w, \tau' = N\tau + N - 1$ .}
\end{lemma}

\medskip

\noindent{Proof:} Let $\overrightarrow{\alpha} =
(\alpha_1,\dots,\alpha_N, 1), \overrightarrow{\beta} =
(\tilde{\alpha_1},\dots,\tilde{\alpha_N}, 1) \in
\mathbb{R}^{N+1}$. We can write
$$
\overrightarrow{\beta} = \frac{1}{\alpha_w} A \overrightarrow{\alpha}
$$
where $A \in GL(N+1, \mathbb{Z})$ is a matrix whose largest entry
is bounded by $1/\alpha_w$, hence $||A^T|| \ll 1 / \alpha_w$ (where
here $||\cdot||$ denotes the operator norm).
Let $\tau' = N\tau + N - 1$. For $\overrightarrow{k} = (p_1,\dots,p_N,q)$
we have
\begin{align*}
|\overrightarrow{\beta} \cdot \overrightarrow{k}| & =
\frac{1}{\alpha_w}|\overrightarrow{\alpha} \cdot A^T
\overrightarrow{k}| \\
& \gg \frac{C}{\alpha_w} \frac{1}{||A^T
\overrightarrow{k}||^{\tau'}} \\
& \geq \frac{C}{\alpha_w} \frac{1}{||A^T||^{\tau'}}
\frac{1}{||\overrightarrow{k}||^{\tau'}} \\
& \gg \frac{C {\alpha_w}^{\tau' - 1}} {||\overrightarrow{k}||^{\tau'}} \\
\end{align*}
The conclusion of the lemma follows from the remarks above.
$\diamond$

\medskip

\begin{lemma}{For $(\alpha_1,\dots,\alpha_N) \in DC_N(C, \tau)$,
given $w \in \{1,\dots,N\}$ and $(\tilde{\alpha_1}, \dots,
\tilde{\alpha_N}) = G(\alpha_1,\dots,\alpha_N, w)$, there exists
$1 \leq j = j(w) \leq N$ such that $\tilde{\alpha_j} \gg C
\alpha^{\tau}_w$.}
\end{lemma}

\medskip

\noindent{Proof:} Let $1/\alpha_w = q + \beta$ where $q \geq 1$ is
an integer and $|\beta| = \tilde{\alpha_w} < 1/2$. We note that
$k_1 \leq q/(q + \beta) \leq k_2$ for universal constants $0 < k_1 < k_2 < 1$.
If $\tilde{\alpha_w} \alpha_w = |1 - q\alpha_w| \geq C/|q|^{\tau}$
then $\tilde{\alpha_w} \geq C/|q|^{\tau}$ and we put $j(w) = w$.
Otherwise we must have $|q\alpha_k - a_k| \geq C/|q|^{\tau}$ for
some $k \neq w$, so
\begin{align*}
\tilde{\alpha_k} = \frac{|a_k \alpha_w - \alpha_k|}{\alpha_w} =
|a_k - \alpha_k(q+\beta)| & \geq |q \alpha_k - a_k| - \alpha_k
|\beta| \\
& \geq \frac{C}{|q|^{\tau}} - \alpha_k |\beta| \\
& \geq {C'}{|q+\beta|^{\tau}} - \alpha_k |\beta| \\
\end{align*}
where $C' \gg C$. If $\alpha_k |\beta| \leq \frac{1}{2}
C'/|q+\beta|^{\tau}$ then $\tilde{\alpha_k} \geq \frac{1}{2} C' \alpha^{\tau}_w$ and
we put $j(w) = k$; otherwise we have
$$
\tilde{\alpha_w} = |\beta| > \frac{1}{\alpha_k} \frac{1}{2} C'
\alpha^{\tau}_w > \frac{1}{2} C' \alpha^{\tau}_w
$$
and we put $j(w) = w$. $\diamond$

\medskip

\begin{prop}{If \ $(\alpha_1,\dots,\alpha_N) \in DC_N(C,\tau)$ \ for some \ $C, \tau >
0$ \ then \ $\mathcal{B}(\alpha_1,\dots,\alpha_n) < +\infty$.}
\end{prop}

\medskip

\noindent{Proof:} We choose a $w \in \{1,\dots,N\}$ and define
$(\alpha^{(n)}_1,\dots,\alpha^{(n)}_N)_{n \geq 0}, (w(n))_{n \geq 0}, (C_n)_{n \geq
0}$ inductively by $(\alpha^{(0)}_1,\dots,\alpha^{(0)}_N) := (\alpha_1,\dots, \alpha_N), w(0) =
w, C_0 = C$, and for $n \geq 0$, $(\alpha^{(n+1)}_1,\dots,\alpha^{(n+1)}_N) :=
G(\alpha^{(n)}_1,\dots,\alpha^{(n)}_N, w(n))$; by Lemma \ref{dioph},
$(\alpha^{(n+1)}_1,\dots,\alpha^{(n+1)}_N) \in DC_N(C_{n+1}, \tau)$, where
$C_{n+1} := k C_n \alpha^{(n) (\tau' - 1)}_{w(n)}$ for some constant $k > 0$ and
$\tau' = N\tau + N - 1$. We let
$w(n+1) := j(w(n))$ be given by the previous lemma applied to
$(\alpha^{(n)}_1,\dots,\alpha^{(n}_N)$ and $w(n)$, so that (for
$k$ taken small enough) $\alpha^{(n+1)}_{w(n+1)} \geq k C_n \alpha^{(n)
\tau}_{w(n)}$.

\medskip

An easy induction gives $C_{n+1} = k^{n+1} C (\alpha^{(0)}_{w(0)}
\dots \alpha^{(n)}_{w(n)})^{\tau' - 1}$. Putting $\tau'' = \max(\tau, \tau' - 1)$
it follows that
$$
\frac{1}{\alpha^{(n+1)}_{w(n+1)}} \leq \frac{1}{k^{n+1} C} \frac{1}{\left(\alpha^{(0)}_{w(0)} \dots
\alpha^{(n)}_{w(n)}\right)^{\tau''}}
$$
so
\begin{align*}
\sum_{n = 0}^{\infty} \alpha^{(0)}_{w(0)} \dots \alpha^{(n)}_{w(n)} \log
\frac{1}{\alpha^{(n+1)}_{w(n+1)}} & \leq \sum_{n = 0}^{\infty} \left(\alpha^{(0)}_{w(0)} \dots
\alpha^{(n)}_{w(n)}\right) \left((n+1) \log \frac{1}{k} + \log \frac{1}{C} + \tau''
\log \frac {1}{\left(\alpha^{(0)}_{w(0)} \dots
\alpha^{(n)}_{w(n)}\right)} \right) \\
& < +\infty
\end{align*}
since $\alpha^{(k)}_{w(k)} < 1/2$ for $k \geq 1$. $\diamond$


\bibliography{simultlin}
\bibliographystyle{alpha}

\medskip

\noindent Ramakrishna Mission Vivekananda University,
Belur Math, WB-711202, India

\noindent email: kingshook@rkmvu.ac.in

\end{document}